\documentclass{amsart}
\usepackage{amssymb}
\usepackage{amsmath}
\usepackage[american]{babel}

\makeatletter
\@addtoreset{equation}{section}
\makeatother

\renewcommand\thefigure{\thesection.\@arabic\c@figure}
\renewcommand\thetable{\thesection.\@arabic\c@table}

\newtheorem{theorem}{Theorem}[section]
\newtheorem{lemma}[theorem]{Lemma}
\newtheorem{proposition}[theorem]{Proposition}

 \theoremstyle{remark}

\newcommand{\mc}[1]{{\mathcal #1}}

\newcommand{\bb}[1]{{\mathbb #1}}

\newcommand{\<}{\langle}
\renewcommand{\>}{\rangle}

\DeclareMathOperator{\spa}{span}

\title[Fluctuations on the Sierpinski gasket]{Equilibrium fluctuations for the zero-range process on the Sierpinski gasket}

\author{M.D. Jara}

\address{\noindent IMPA, Estrada
Dona Castorina 110, CEP 22460 Rio de Janeiro, RJ, Brazil
\newline
e-mail:  \rm \texttt{monets@impa.br}}

\date{}

\begin{document}

\begin{abstract}
We obtain the fluctuations from the hydrodynamic limit for the zero-range process in the Sierpinski gasket $V$. The limiting process is given by a generalized Ornstein-Uhlenbeck process associated to the Neumann Laplacian on $V$.

\end{abstract}

\maketitle

\section{Introduction}

A considerable amount of progress has been accomplished in the study of the motion of a particle in fractal structures, like the Sierpinski gasket or the Sierpinski carpet. See \cite{Kig}, \cite{Bar} for a survey of the field. However, much less is known for systems of interacting particles evolving in fractal structures. In the other hand, a well-developed theory of hydrodynamic limit of interacting particle systems has been developed \cite{KipLan}.  

In this article, we perform a first step into the study of interacting particle systems in fractal structures. We obtain the fluctuations of the empirical distribution of particles for the zero-range process from its equilibrium state. We prove that, when properly rescaled, the fluctuations are given by a generalized Ornstein-Uhlenbeck process defined on a Sobolev-like space of functionals over the Sierpinski gasket. In particular, the fluctuations are given by a Gaussian process, like in the case of the cubic lattice, but are sub-diffusive in the sense that the scaling exponent is not 2. This collective behavior of the particles is different from the behavior of a single particle. In fact, the scaling limit of a  simple random walk on the graph approximations to this fractal is the so-called Brownian motion on the Sierpinski gasket \cite{BarPer}, and this process is not a semimartingale.

We organize this article as follows. In section \ref{s1}, we introduce the Laplacian on the Sierpinski gasket and we construct the generalized Ornstein-Uhlenbeck process as an appropriated limit of finite-dimensional Ornstein-Uhlenbeck processes. Notice that this construction applies each time the Laplacian has a compact resolvent. In particular, the results of this article can be obtained in a straight-forward way for any finitely ramified fractal in the sense of Kigami \cite{Kig}. In section \ref{s2} we introduce the zero-range process and we outline the proof of the main result. In section \ref{s3} we complete the proof of Theorem \ref{t2} obtaining the necessary lemmas.

\section{Generalized Ornstein-Uhlenbeck process in the Sierpinski gasket}
\label{s1}
\subsection{The Sierpinski gasket}
Let $a_0 =(0,0)$, $a_1=(1,0)$, $a_2=(1/2,\sqrt 3/2)$ be the vertexes of the equilateral triangle $V_0=\{a_0,a_1,a_2\}$. Define $f_i(z)= (z+a_i)/2$ for $i=0,1,2$. The Sierpinski gasket is the unique compact, non-empty set $V \subseteq \bb R^2$ such that 
\[
 V = \bigcup_{i=0,1,2} f_i(V).
\]

A more explicit, recursive construction of $V$ is the following. 
For each $n \geq 0$, define inductively $V_{n+1} = \cup_{i=0,1,2}f(V_n)$, and set $V^*=\cup_{n \geq 0} V_n$. 
Then, $V$ is the closure of $V^*$ in $\bb R^2$. Let $\mu_n$ denote the measure which assigns mass $3^{-n}$ to each point on $V_n$. 
The sequence of measures $\{\mu_n\}$ converge in the vague topology to the Hausdorff measure $\mu$ of dimension $d_f=\log 3/\log 2$ on $V$.
We define a graph $\Gamma_n=(V_n,E_n)$ in the natural way: $e_n=\{\<xy\>; x, y \in V_n, |x-y| =2^{-n}\}$. 
For each function $f: V^* \bb R$ and each $n \geq 0$ we define 
\[
 \mc E_n(f,f) = \big(5/3\big)^n \sum_{\<xy\> \in E_n} \big(f(y)-f(x)\big)^2.
\]
A simple computation shows that for any $f$ the sequence $\{\mc E_n(f,f)\}_n$ is non-decreasing. 
Therefore, the (possibly infinite) limit
\[
 \mc E(f,f) = \lim_{n \to \infty} \mc E_n(f,f)
\]
is well-defined. We say that $f \in \mc H_1$ if $\mc E(f,f) < \infty$. For any function $f: V^* \to \bb R$, we have the folowing estimate \cite{Kig}:
\begin{equation}
\label{ec3}
 \sup_{x,y \in V_n} \frac{|f(x)-f(y)|}{|x-y|^\alpha} \leq 6 \sqrt{\mc E_n(f,f)},
\end{equation}
where $\alpha = (\log 5/3)/\log4$.
Therefore, any function $f \in \mc H_1$ is H\"older-continuous and can be continuously extended to $V$. For each $f$ in $\mc H_1$, we let $\bb H_n f$ be the solution of the variational problem
\[
 \inf\{\mc E(g,g); g(x) =f(x) \text{ for all } x \in V_n\}.
\]

A simple computation shows that $\bb H_n f$ is well-defined. In particular, the space $\mc H_1$ is non-empty and dense in $\mc C(V)$, the set of continuous functions $f: V \to \bb R$. We give a Hilbert structure to $\mc H_1$ by taking the norm $||f||_1= \mc E(f,f)^{1/2}$ and defining the inner product $\mc E(f,g)$ by polarization. We have the inclusions $\mc H_1 \subseteq \mc C(V) \subseteq \mc L_2(\mu)$. Define $\mc H_{-1}$ as the dual space of $\mc H_1$ with respect to $\mc L_2(\mu)$: for $f \in \mc L_2(\mu)$, 
\[
 ||f||_{-1} := \sup\Big\{\int fg d\mu; g \in \mc H_1, ||g||_1=1\Big\}
\]
and let $\mc H_{-1}$ be the closure of $\mc L^2(\mu)$ with respect to this norm.
Notice that for any $f,g \in \mc H_1$, $ \mc E_n(f,g) = -\int f \Delta_n g d\mu_n$, where $\Delta_n$ is the operator given by
\[
 \Delta_n g(x) = 5^n \sum_{y:\<xy\> \in E_n} \big(f(y)-f(x)\big).
\]
The operator $\Delta_n$ corresponds to the discrete Laplacian on the graph $\Gamma_n$ with some appropriated weight. This fact motivates the definition of the Laplacian in $V$ as the generator of the inner product $||f||_1$ with respect to the measure $\mu$: for $g \in \mc H_1$, we define $\Delta g$ as the unique element of $\mc H_{-1}$ such that $\int f \Delta g d\mu= \mc E(f,g)$ for all $f \in \mc H_1$. We can think indistinctly the operator $\Delta$ as an isomorphism from $\mc H_1$ to $\mc H_{-1}$ or as an unbounded, densely defined operator from $D(\Delta) \subseteq \mc L_2(\mu)$ to $\mc L_2(\mu)$. 
The operator $\Delta$ is called the {\em Neumann Laplacian} on $V$. In fact, the operator $-\Delta$ is essentially self-adjoint, non-negative and has a discrete spectrum $\{\lambda_k, k \in \bb N_0\}$, with $0= \lambda_0 < \lambda_1\leq ... < +\infty$ and $\lim_k \lambda_k = +\infty$. 
An extensive treatment of these results based on explicit computations for the case of the Dirichlet Laplacian can be found in \cite{FS}. A proof of these results for the more general class of p.c.f.s.s sets can be found in \cite{Kig}.

\subsection{Construction of the generalized Ornstein-Uhlenbeck process}
Let $v_k$ be a normal eigenvector associated to the eigenvalue $\lambda_k$. Then, $\{v_k, k\geq 0\}$ is an orthonormal basis for $\mc L_2(\mu)$ and also for $\mc H_1$, $\mc H_{-1}$. 
For a function $f \in \mc L_2(\mu)$, define the Fourier coefficients $f_k$ by $f_k = \int fv_k d\mu$. We define the Sobolev norm $||\cdot||_m$, $m \in \bb R$ by
\[
 ||f||_m^2 = \sum_{k \geq 0} \lambda_k^m f_k^2.
\]

We denote by $\mc H_m$ the completion of $\{f; ||f||_m < \infty\}$ under this norm. Notice that this definition is consistent with our previous definition of $\mc H_1$, $\mc H_{-1}$, and $\mc H_0 = \mc L_2(\mu)$. Now we can state the main result of this section.

\begin{theorem}
\label{t1}
Fix $\beta, \gamma > 0$ and $m \geq 2$. There exists a process $\mc Y_t$ with trajectories in the space of continuous paths $\mc C([0,\infty), \mc H_{-m})$ that solves the following martingale problem:
\begin{itemize}
\item[{\bf (MP)}] 
For any function $g \in \mc H_m$. 
\[
 \mc M_t(g) = \mc Y_t(g) -\mc Y_0(g) - \int_0^t \beta \mc Y_s(\Delta g) ds
\]
is a martingale of quadratic variation
\[
 \<\mc M_t(g)\> = \gamma t\mc E(g,g)
\]
 
\end{itemize}

Moreover, $\mc Y_t$ is uniquely determined by its initial  distribution.
\end{theorem}

Define $H_k = \spa \{v_0,...,v_k\}$. Clearly, $H_k$ is isomorphic to $\bb R^{k+1}$. We define $\mc Y_t^{k,0}$ as the It\^o's process in $\bb R^{k+1}$, solution of the stochastic equation
\[
 d \mc Y_i^{k,0}(t) = -\beta \lambda_i \mc Y_i^{k,0}(t) dt + \sqrt{\gamma \lambda_i} d B_i(t), i=0,...,k,
\]
\[
 \mc Y_i^{k,0}(0) = Y_i^{k,0},
\]
where $(B_0(t),...,B_k(t))$ is a $k+1$-dimensional standard Brownian motion and $Y^{k,0}=(Y_0^{k,0},...,Y_k^{k,0})$ is the initial condition. We define then the process $\mc Y_t^k$ in by
\[
 \mc Y_t^k(f) = \sum_{i=0}^k f_i \mc Y_i^{k,0}(t),
\]

Notice that $\mc Y_t^k$ is well defined in $\mc H_{m}$ for every $m \in \bb R$.
If the sequence of random variables $Y_0^k(0)$ is consistent, then the sequence of processes $\mc Y_i^{n,0}(t)$ is also consistent, and by Kolmogorov's extension theorem there exists a limiting process for this sequence. However, this limiting process does not have {\em a priori} continuous trajectories. We will prove that the sequence of processes $\{ \mc Y_t^k\}_k$ is tight in $\mc C([0,T],\mc H_{-m})$ for every $T >0$. By Mitoma's criterion, it is enough to prove that $\{\mc Y_t^k(f)\}_k$ is tight in $\mc C([0,T],\bb R)$ for all $f \in \mc H_m$, and that for every $\epsilon >0$,
\begin{equation}
\label{ec4}
\lim_{\delta \to 0} \sup_{f: ||f||_m \leq \delta} P\Big[\sup_{t \in [0,T]} \big|\mc Y_t^k(f)\big| \geq \epsilon\Big] = 0.
\end{equation}

By It\^o's formula,
\begin{equation}
\label{ec1}
 \mc M_t^k(f) = \mc Y_t^k(f) - \mc Y_0^k(f) - \beta\int_0^t \mc Y_s^k(\Delta f) ds
\end{equation}
 is a martingale of quadratic variation 
 \[
  \<\mc M_t^k(f)\> = \gamma t \sum_{i=0}^k \lambda_i f_i^2.
 \]
In particular, $\mc M_t^k$ is a Brownian motion of mean zero and variance $\gamma \sum_{i=1}^k \lambda_i f_i^2$. Therefore, the sequence $\{\mc M_t^k(f)\}_k$ is tight if and only if 
\[
 \sup_k \sum_{i=1}^k \lambda_i f_i^2 = \sum_{i=0}^\infty \lambda_i f_i^2 < +\infty.
\]
Now we show that the integral term in \ref{ec1} is tight. It is enough to prove that 
\[
 \sup_k \sup_{t \in [0,T]} E[\mc Y_t^k(\Delta f)^2] < \infty.
\]
Define $\psi_i(t)= E[\mc Y_i^{k,0}(t)^2]$, and assume that $\psi_i(0) < \infty$. Applying It\^o's formula to $\mc Y_i^{k,0}(t)^2$, we obtain that
\[
 \psi_i(t) = \frac{\gamma}{2\beta}\Big( 1 -e^{-2\beta \lambda_i t}\Big) + \psi_i(0) e^{-2\beta \lambda_i t}.
\]
In the same way, defining $\psi_{ij}(t) = E[\mc Y_i^{k,0}(t)\mc Y_j^{k,0}(t)]$, we see that
\[
 \psi_{ij}(t) = \psi_{ij}(0) e^{-\beta(\lambda_i +\lambda_j)t}.
\]

By the definition of $\mc Y_t^k(f)$, we see that
\[
 E[\mc Y_t^k(\Delta f)^2] = \sum_{i=0}^k  \frac{\gamma\lambda_i^2 f_i}{2\beta} \Big(1- e^{-2\beta \lambda_i t}\Big) +(\Delta P_t f^k)^* \Psi_k (\Delta P_t f^k), 
\]
where $\Delta P_t f^k$ is an abbreviation for the vector $(\lambda_0 e^{-\beta \lambda_0 t} f_0,...,\lambda_k e^{-\beta \lambda_k t} f_k)$, and $\Psi_k$ is the matrix $(E[\mc Y_i^{k,0}(0)\mc Y_j^{k,0}(0)])_{ij}$. In particular, the integral term in \ref{ec1} is tight, provided that
\[
 \sum_{i \geq 0} \lambda_i^2 f_i^2 < +\infty \text{ and } \sup_{k \geq 0} ||\Psi_k|| < +\infty.
\]

Since $\lambda_1>0$, this last condition imply that $\sum_i \lambda_i f_i^2 <+\infty$. This is equivalent to $f \in \mc H_2$. Let $\mc Y_0$ be the limiting distribution of $\mc Y_0^k$. The condition in $\Psi$ is equivalent to the condition $E[\mc Y_0(f)^2] \leq C||f||_2$ for every $f \in \mc H_2$. Using Tchebyshev's and Doob's inequalities, we see that these very same conditions allow us to obtain the estimate \ref{ec4}.
Therefore, we have proved the following result.

\begin{theorem}
 Let $\mc Y_0$ be a random variable in $\mc H_{-2}$ such that there exists a constant $C$ with $E[\mc Y_0(f)^2] \leq C||f||_2$ for every $f \in \mc H_2$. Then, the processes $\mc Y_t^k$ in $\mc C([0,T],\mc H_{-2})$ form a tight sequence.
\end{theorem}

Now we want to prove that the limit point is unique. Let $\mc Y_t$ be a limit point of $\mc Y_t^k$. For each $f \in H_2$, taking the approximation $f_{(n)}= \sum_{i \leq n} f_i v_i$, it is easy to prove that 
\[
 \mc M_t(f) = \mc Y_t(f) - \mc Y_0(f) - \int_0^t \beta \mc Y_s(\Delta f) ds,
\]
\[
 \big(\mc M_t(f)\big)^2 - \gamma t ||f||_1^2
\]
are martingales (this is straightforward if $f_n=0$ for $n$ large enough). Take $t_1,...,t_n \in [0,T]$ and $f^{(1)},...,f^{(n)}$ in $\mc H_2$. Using this martingale representation, it is easy to obtain the distribution of the vector $(\mc Y_{t_1}(f^{(1)}),...,\mc Y_{t_n}(f^{(n)}))$. Since these finite-dimensional distributions determine probability in $\mc C([0,T], \mc H_{-2})$, we conclude that the limit point $\mc Y_t$ is unique, and satisfies the martingale problem of Theorem \ref{t1}.

\section{The zero-range process in the Sierpinski gasket}
\label{s2}
Consider the state space $\Omega_n = \bb N_0^{V_n}$. For any element $\eta \in \Omega_n$ and any site $x \in V_n$, $\eta(x)$ represents the number of particles at $x$. Let $h: \bb N_0 \to \bb R_+$ be a given function with $h(0)=0$. Particles evolve in $V_n$ according to the following rule. At each time $t$, a particle jumps from site $x$ to its neighbor $y$ at an exponential rate $h(\eta_t(x))$. This happens independently for each pair $\<xy\>$ of neighboring sites. After a jump, a new exponential clock begins with the corresponding rate, independently of the past of the process. In this way we have described a Markov process $\eta_t$ with state space $\Omega_n$ an generated by the operator
\[
 \mc L_n f(\eta) = \sum_{\<xy\> \in E_n} h\big(\eta(x)\big) \big[ f(\eta^{xy}) -f(\eta)\big],
\]
where 
\[
 \eta^{xy}(z) = \left\{ \begin{array}{cl}
                 \eta(x)-1, &z=x\\
		 \eta(y)+1, &z=y\\
		 \eta(z), &z \neq x,y.\\
                \end{array}
		\right.
\]

We will assume that the jump rate $h(\cdot)$ has linear growth and bounded variation: there exists a constant $\epsilon_0>0$ such that $\epsilon_0 n \leq h(n) \leq \epsilon_0^{-1}$, and $\sup_n|h(n+1) - h(n)| < +\infty$. With these assumptions, the process $\eta_t$ has  a one-parameter family of invariant measures $\nu_\rho$, $\rho \in [0,\infty)$ of product form. The marginals of $\nu_\rho$ are given by
\[
 \nu_\rho(\eta(x) = k) = \frac{1}{Z(\rho)}\frac{\phi(\rho)^k}{h(k)!},
\]
where $h(k)!= h(1)\cdots h(k)$, $h(0)!=1$, $Z(\rho)$ is a normalization constant and $\phi(\rho)$ is chosen in such a way 
that $\int \eta(x) d\nu_\rho(\eta)=\rho$ for every $x \in V_n$. Notice that $\int h\big(\eta(x)\big) d\nu_\rho= \phi(\rho)$. 

\subsection{The fluctuation field}
For each function $f: V^* \to \bb R$, consider the empirical measure $\pi_t^n(f)$  by
\[
 \pi_t^n(f) = \frac{1}{3^n} \sum_{x \in V_n} \eta_t(x) f(x),
\]
where $\eta_t$ is the zero-range process defined in the previous section. When the process $\eta_t$ starts from the equilibrium measure $\nu_\rho$, by the law of large numbers $\pi_t^n(f)$ converges to $\rho \int f d\mu$ for any continuous function $f: V \to \bb R$.  
In order to investigate the fluctuations of the empirical measure, we define the fluctuation field $\mc Z_t^n$ by its action over functions $f$ in some appropriated Sobolev space $\mc H_{m}$:
\[
 \mc Z_t^n(f) = \frac{1}{3^{n/2}} \sum_{x \in V_n} \big(\eta_t^n(x) - \rho\big) f(x),
\]
where $\eta_t^n= \eta_{5^nt}$. We have speeded up the process by $5^n$ in order to obtain a non-trivial limit when $n \to \infty$. 

\begin{theorem}
 \label{t2}
 Fix a positive integer $m \geq 3$. Then the process  $\mc Z_t^n$ converges in distribution to the Ornstein-Uhlenbeck process $\mc Y_t$ with characteristics $\beta = \phi'(\rho)$ and $\gamma = \phi(\rho)$ and initial distribution $\mc Y_0$, where $\mc Y_0$ is a Gaussian field of mean zero and covariance $E[\mc Y_0(f) \mc Y_0(g)]= \int fg d\mu$.
\end{theorem}

Now we outline the proof Theorem \ref{t2}, stating the main steps as lemmas to be proved in the next sections.
\begin{proof}
For each function $f$, Dynkin's formula shows that
\begin{equation}
 \label{ec2}
 M_t^n(f) = \mc Z_t^n(f)- \mc Z_0^n(f) - \int_0^t \frac{1}{3^{n/2}} \big(h\big(\eta_s^n(x)\big)- \phi(\rho)\big) \Delta_n f(x) ds
\end{equation}
is a martingale of quadratic variation 
\[
 \<M_t^n(f)\> = \int_0^t \frac{5^n}{3^n} \sum_{\<xy\> \in E_n} \big(h\big(\eta_s^n(x)\big) + h\big(\eta_s^n(y)\big)\big)\big[ f(y)- f(x)\big]^2 ds.
\]

For any $m \geq 3$, by Lemma \ref{l1} the sequence of processes $\{ \mc Z_t^n\}_n$ is tight in $\mc D([0,T],\mc H_{-m})$, and the limit points are concentrated in continuous trajectories. Let $\mc Z_t$ a limit point of this sequence. To simplify the notation, we just call $n$ the subsequence for which $\mc Z_t^n$ converges to $\mc Z_t$. By Lemma \ref{l3}, we can replace  $h\big(\eta_s^n(x)\big)-\phi(\rho)$ by $\phi'(\rho)(\eta_s^n(x)-\rho)$ in \ref{ec2}. By Lemma \ref{l2}, for any function $f \in \mc H_m$ the martingales $\mc M_t^n(f)$ converge in distribution to the martingale
\[
 \mc Z_t(f) - \mc Z_0(f) - \int_0^t \phi'(\rho) \mc Z_s(\Delta f) ds
\]
of quadratic variation $\phi(\rho)t \mc E(f,f)$. By the central limit theorem and standard computations, the initial distribution $\mc Z_0^n$ converges to $\mc Y_0$. Therefore, the process $\mc Z_t$ satisfy the martingale problem \ref{ec1}. Since there is a unique process satisfying this martingale problem, we conclude that $\mc Z_t= \mc Y_t$. Therefore, the sequence $\mc Z_t^n$ has a unique limit point. Since the topology of convergence in distribution is metrizable, we conclude that the whole sequence $\mc Z_t^n$ converges in distribution to $\mc Y_t$.

\end{proof}

\section{Proofs}
\label{s3}
In this section we prove Lemmas \ref{l1}, \ref{l2}, \ref{l3}, closing in this way the proof of Theorem \ref{t2}.
\subsection{Tightness of $\mc Z_t^n$}
In this subsection we prove tightness of the sequence $\{\mc Z_t^n\}_n$ in $\mc D([0,T],\mc H_{-m})$. First, notice that $\mc Z_t^n$ is in $\mc H_{-1}$. In fact, by equation \ref{ec3}, for any function $f \in \mc H_1$, $||f||_\infty \leq ||f||_0+6||f||_1$, from which the Dirac-$\delta$ distribution, and therefore $\mc Z_t^n$ is in $\mc H_{-1}$. By Mitoma's criterion, we have to prove that $\mc Z_t^n(f)$ is tight for any $f \in \mc H_m$, and that for any $\epsilon >0$, 
\[
\lim_{\delta \to 0} \sup_{||f||_m \leq \delta} \sup_{n \geq 0} \bb P_\rho\big( \sup_{t \in [0,T]} \big| \mc Z_t^n(f)\big| \geq \epsilon\big) = 0.
\]

By the martingale decomposition of $\mc Z_t^n(f)$, it is enough to prove the following estimates.
\begin{proposition}
 The sequence $\mc Z_t^n$ of processes with trajectories in $\mc D([0,T], \mc H_{-m})$ is tight if for any $f$ in $\mc H_m$, we have
 \begin{itemize}
\item[i)]
\[
 \sup_{n \geq 0} \sup_{t \in [0,T]} \bb E_\rho\Big[ \Big(\frac{1}{3^{n/2}} \sum_{x \in V_n} \big(h(\eta_t^n(x)) -\rho\big) \Delta_n f(x)\Big)^2\Big]=0
\]
\item[ii)] \[
             \sup_{n \geq 0} \sup_{t \in [0,T]} \bb E_\rho\Big[ \Big(\frac{5^n}{3^n} \sum_{\<xy\>  \in E_n} h(\eta_t^n(x)) \big(f(y)-f(x)\big)^2\Big)^2\Big]=0. 
            \]
 \end{itemize}
Moreover, the limit points of $\mc Z_t^n$ are concentrated on continuous trajectories if
\begin{itemize}
\item[iii)] There exists $\delta(f,n)$ such that $\lim_{n\to \infty} \delta(f,n)=0$ and 
\[
	\lim_{n \to \infty} \bb P_\rho \Big[\sup_{t \in [0,T]} \big|\mc Z_{t+}^n(f)-\mc Z_t^n(f)\big| > \delta(f,n)\Big]=0.
\]
\end{itemize}
\end{proposition}

Using the invariance and the product form of the initial distribution $\nu_\rho$, it is not hard to rewrite conditions $i)$ and $ii)$ in terms of conditions on $f$. 
The first estimate follows if $\Delta f$ is continuous. Since functions in $\mc H_1$ are continuous, for $m \geq 3$, $\Delta f$ is continuous. This is the only point where we need $m$ to be larger than 3. The second estimate follows when $f \in \mc H_1$ from the fact that functions $f$ in $\mc H_1$ are ``differentiable'' in a classical variational sense:
\[
\lim_{n \to \infty} \sup_{x \in V_n} \big(5/3\big)^n \sum_{y:\<xy\> \in E_n} \big(f(y)-f(x)\big)^2 = 0.
\]

Since the process $\eta_t^n$ has no simultaneous jumps for almost every trajectory, the fourth estimate follows taking $\delta(f,n) = 2||f||_\infty/ 3^{n/2}$. Therefore, we have proved the following result:

\begin{lemma}
 \label{l1}
 For any $m \geq 3$, the sequence $\{\mc Z_t^n\}_n$ of processes in $\mc D([0,T], \mc H_{-m})$ is tight. 
\end{lemma}

\subsection{The Boltzmann-Gibbs principle}
In this subsection we prove the following lemma, known as the Boltzmann-Gibbs principle:
\begin{lemma}[Boltzmann-Gibbs principle]
For any continuous function $f: V \to \bb R$, 
\[
\lim_{n \to \infty} \bb E_\rho\Big[ \Big(\int_0^t \frac{1}{3^{n/2}} \sum_{x \in V_n} \big[h(\eta_s^n(x)) - \phi(\rho) - \phi'(\rho) \big(\eta_s^n(x) - \rho\big)\big] f(x) ds\Big)^2\Big] =0.
\]

\label{l2}
 
\end{lemma}
What this theorem says is that the fluctuations of non-conserved quantities are much faster than the fluctuations of  conserved quantities. Therefore, when properly rescaled, the only component of the non-linear fluctuation field associated to $h(\eta(x))$ that survives is the projection over the conserved quantity, in this case, the number of particles $\eta(x)$. The Boltzmann-Gibbs principle was first introduced by Rost \cite{Ros}. We adopt here Chang's proof \cite{Cha} (see also \cite{KipLan}).

\begin{proof}
The fist step into this proof is to localize the problem. Let $k>0$ be a fixed integer, that will go to $\infty$ after $n$. For each point $x \in V_n$, we define $\Delta_n^k(x) \subseteq V_n$ as the triangle of side $2^{-(n-k)}$ that contains $x$ and with vertexes in $V_{n-k}$. Notice that for $x \in V_{n-k} \setminus V_0$, there are two possible choices for $\Delta_n^k(x)$, and exactly one choice for $x \notin V_{n-k}\setminus V_0$. We solve this tie in any consistent way, it does not really matter, after all. Let $b_k = 3(3^k+1)/2$ be the number of points from $V_n$ in $\Delta_n^k(x)$. Notice that $b_k$ does not depend on $n$. Define $\mc V_{n,k}(x,\eta)$ by
\[
\mc V_{n,k}(x,\eta) = \frac{1}{b_k} \sum_{y \in \Delta_n^k(x)} \Big(h(\eta(y)) - \phi(\rho) -\phi'(\rho)(\eta(y)-\rho)\Big).
\]

There are $3^{n-k}$ possible choices for the sets $\Delta_n^k(x)$. For each one of them, choice a representant $x_i^n$. The sequence $\{x_i^n\}_i$ is chosen in such a way that for any $x \in V_n$ there exists $i$ with $\Delta_n^k(x) = \Delta_n^k(x_i^n)$, and for any $i\neq j$, $\Delta_n^k(x_i^n) \neq \Delta_n^k(x_j^n)$. Since $f$ is continuous and $V$ is compact, $f$ is uniformly continuous. Therefore, the theorem follows if we prove that 
\begin{equation}
\label{ec5}
\lim_{k \to \infty} \lim_{n \to \infty} \bb E_\rho \Big[\Big(\int_0^t \frac{1}{3^{n/2}} \sum_{i=1}^{3^{n-k}} \mc V_{n,k}(x_i^n, \eta_s^n) f(x_i^n)ds\Big)^2\Big]=0.
\end{equation}

Here we have introduced a convenient averaging on $\Delta_n^k(x_i)$, that we will exploit below. Let $\mc L_{n,i}$ the generator of the process restricted to $\Delta_n^k(x_i)$. Take a function $\Theta: \bb N_0^{\Delta_n^k(x_1)}$ in $\mc L_2(\nu_\rho)$, and let $\Theta_i$ be its translation to $\bb N_0^{\Delta_n^k(x_i)}$. Notice that the dynamics on each block $\Delta_n^k(x_i)$ is the same. Using the estimate \cite{KipLan}
\[
\bb E_\rho\Big[ \Big(\int_0^t \varphi(\eta_s^n) ds\Big)^2\Big] 
	\leq 20 t \int \varphi(\eta)(-5^n\mc L_n)^{-1} \varphi(\eta) d\nu_\rho,
\]
we obtain that
\[
\bb E_\rho\Big[\Big(\int_0^t \frac{1}{3^{n/2}} \sum_{i=1}^{3^{n-k}} \mc L_{n,i} \Theta_i(\eta_s^n) f(x_i^n) ds\Big)^2\Big] 
	\leq \frac{20t||f||^2_\infty}{3^k5^n} \int \Theta(\eta) (-\mc L_{n,0})\Theta(\eta) d\nu_\rho,
\]

where we have also used the fact that the dynamics is the same on each block $\Delta_n^k(x_i)$.
In particular, the left-hand side of this inequality converges to $0$ as $N \to \infty$. Therefore, we can substract this term from equation \ref{ec5}, and it is enough to prove that

\begin{multline*}
\lim_{k \to \infty} \inf_{\Theta \in \mc L_2(\nu_\rho)} \limsup_{n \geq 0}\\
\bb E_\rho\Big[\Big(\int_0^t \frac{1}{3^{n/2}} \sum_{i=1}^{3^{n-k}} f(x_i^n)\big[\mc V_{n,k}(x_i^n,\eta_s^n) - \mc L_{n,i} \Theta_i(\eta_s^n)\big]ds\Big)^2\Big]=0.
\end{multline*}

Here the infimum is over functions defined in $\bb N_0^{\Delta_n^k(x_0^n)}$. Remember that $\nu_\rho$ is a product measure. Therefore, integrals over non-adjacent triangles are independent. Putting  the expectation in the previous integral into the time integral, using stationarity and letting $N$ go to $\infty$, we see that it is enough to prove that
\[
\lim_{k \to \infty} \inf_{\Theta \in \mc L_2(\nu_\rho)} \frac{1}{3^k} \int \big( \mc V_k(\eta) - \mc L_k \Theta\big)^2 d\nu_\rho = 0,
\]
where $\mc V_k(\eta)=\mc V_{k,k}(x,\eta)$ and $\mc L_k$ is the generator of the process on the triangle $\Delta_k$. Since $\Delta_k$ has a number of points of order $3^k$, this last limit follows from the equivalence of ensembles and standard arguments. 
\end{proof}
\subsection{The martingale problem}
Once we have proved that the sequence of processes $\{\mc Z_t^n\}_n$ is tight, and we know how to write the integral term in \ref{ec2} as a function of the process $\mc Z_t^n$, we need to prove that the martingale decomposition of $\mc Z_t^n$ still holds  in the limit. Let $\mc Z_t$ be a limit point of $\mc Z_t^n$. We simply denote by $n$ the subsequence for which the convergence is realized. By Lemma \ref{l2}, 
\[
\int_0^t \frac{1}{3^{n/2}} \sum_{x \in V_n} \big( h(\eta_s^n(x)) - \phi(\rho)\big) \Delta_n f(x) ds - \int_0^t \phi'(\rho) \mc Z_s^n(\Delta f) ds \xrightarrow{n \to \infty} 0
\]
in $\mc L_2(\bb P_\rho)$ and, in particular, in probability. Therefore, the integral term in \ref{ec2} converges in distribution to $\int_0^t \phi'(\rho)\mc Z_s(\Delta f) ds$. Since $\mc Z_t^n(f)$ converges in distribution to $\mc Z_t(f)$, we see that $M_t^n(f)$ converges in distribution to 
\[
M_t(f) = \mc Z_t(f) - \phi'(\rho) \int_0^t \mc Z_s(\Delta f) ds.
\]

\begin{lemma}
 \label{l3}
For any $f \in \mc H_m$ with $m \geq 2$, the processes $M_t(f)$ and $M_t(f)^2 - \phi(\rho) T ||f||_1^2$ are martingales.
\end{lemma}
\begin{proof}
A sufficient condition for the limit $M_t(f)$ to be a a martingale is the existence of a deterministic constant $b$ such that the jumps of $M_t^n(f)$ are uniformly bounded by $b$. But the jumps of $M_t^n(f)$ are bounded by $||f||_\infty/3^{n/2}$. Notice now that the quadratic variation $\<M_t^n(f)\>$ converges in $\mc L_2(\bb P_\rho)$ to $\phi(\rho)t||f||_1^2$. Therefore,
\[
M_t^n(f)^2 - \<M_t^n(f)\> \xrightarrow{N \to \infty} M_t(f)^2 - \phi(\rho) t ||f||_1^2
\]
in distribution. Again, since the jumps of $M_t^n(f)$ are uniformly bounded on $n$, we conclude that $M_t(f)^2-\phi(\rho)t||f||_1^2$ is also a martingale.
\end{proof}

\end{document}